\documentclass[a4paper,12pt]{amsproc}
\usepackage{amssymb,amsthm}
\usepackage{graphicx}

\usepackage{cite}
\usepackage{url}

\newcommand{\e}{\varepsilon}
\newcommand{\pa}{\partial}
\newcommand{\be}{\beta}
\newcommand{\ga}{\gamma}
\newcommand{\dt}{\delta}
\newcommand{\la}{\lambda}
\newcommand{\si}{\sigma}
\newcommand{\ka}{\kappa}
\newcommand{\ph}{\varphi}
\newcommand{\wt}{\widetilde}

\newcommand{\Ga}{\Gamma}

\newcommand{\abs}[1]{\lvert#1\rvert}

\newcommand{\pig}{\pi_G^\#}
\newcommand{\Ng}{\mathcal N_G^\#}
\newcommand{\Ngk}[1][k]{\mathcal N_{G,#1}^\#}
\newcommand{\zg}[1][s]{\zeta_G(#1)}
\newcommand{\zgk}[1][s]{\zeta_{G,k}(#1)}

\newcommand{\Sgk}{\sum_{a\in G_k}}

\newcommand{\Pp}{\prod_{p\in P}}

\newcommand{\esa}[1][a]{e^{-\pa(#1)s}}
\newcommand{\esp}{\esa[p]}
\newcommand{\esak}[1][a]{e^{-\pa(#1)s k}}
\newcommand{\espk}{\esak[p]}

\newcommand{\BL}{\biggl}
\newcommand{\BR}{\biggr}
\newtheorem{theorem}{Theorem}
\theoremstyle{definition}
\newtheorem{remark}{Remark}

\begin{document}
\title[Asymptotics of the counting function]%
{Asymptotics of the counting function\\
of $k$-th power-free elements\\
in an arithmetic semigroup}

\author[V.\,L.~Chernyshev]{V.\,L.~Chernyshev$^1$}

\author[D.\,S.~Minenkov]{D.\,S.~Minenkov$^2$}

\author[V.\,E.~Nazaikinskii]{V.\,E.~Nazaikinskii$^{2,3}$}

\thanks{$^1$National Research University Higher School of Economics,
Moscow, Russia}

\thanks{$^2$Ishlinsky Institute for Problems in Mechanics,
Russian Academy of Sciences, Moscow, Russia}

\thanks{$^3$Moscow Institute of Physics and Technology
(State University), Dolgoprudny, Moscow District, Russia}

\thanks{E-mail: \texttt{vchern@gmail.com,
minenkov.ds@gmail.com, nazay@ipmnet.ru}}



\date{}

\begin{abstract}
For any $k\ge2$, we find the asymptotics of the counting function
of $k$-th power-free elements in an additive arithmetic semigroup
with exponential growth of the abstract prime counting function.
This paper continues the authors' earlier research dealing with the
case of $k=\infty$.
\end{abstract}

\maketitle

\section{Introduction and statement of the problem}

Let~$G$ be an additive arithmetical
semigroup~\cite[pp.~11,~56]{Knopf}; i.e.,

(i) $G$ is a~commutative semigroup with identity element~$1$.

(ii) There exists a (uniquely determined) countable subset
$P\subset G$ (whose elements are called the \textit{primes} of $G$)
such that every element $a\in G$, $a\ne1$, has a factorization of
the form
\begin{equation}\label{factor}
    a=p_1^{n_1}p_2^{n_2}\dots p_s^{n_s}
\end{equation}
with some positive integers $s,n_1,\dots,n_s$ and elements
$p_1,\dots,p_s\in P$, and this factorization is unique up to the
order of factors.

(iii) A~mapping $\pa\colon G\to\mathbf{R}$ (called the
\textit{degree mapping}) is given such that $\pa(1)=0$, $\pa(p)>0$
for all $p\in P$, $\pa(ab)=\pa(a)+\pa(b)$ for all $a,b\in G$, and
the number $\Ng(x)$ of elements $a\in G$ such that $\pa(a)\le x$ is
finite for every $x>0$.

For an integer $k\ge2$, an element $a\in G$ is said to be
\textit{$k$-th power-free} if it has no divisors of the form $b^k$,
where $1\ne b\in G$. We denote the set of $k$-th power-free elements
in~$G$ by $G_k\subset G$ and the number of $k$-th power-free
elements of degree $\le x$ by $\Ngk(x)$. It is natural to extend
the definition to $k=\infty$ by setting $G_\infty=G$; then
$\Ngk[\infty](x)=\Ng(x)$.

Our aim is to find the asymptotics of $\Ngk(x)$ as $x\to\infty$
under the assumption that the number $\pig(x)$ of primes $p\in P$
such that $\pa(p)\le x$ has the asymptotics
\begin{equation}\label{e1}
    \pig(x) = \rho x^\gamma e^{x} ( 1 + O(x^{-\dt})), \qquad x
\rightarrow \infty,
\end{equation}
for some $\rho > 0$, $\gamma > -1$, and $\dt\in(0,1]$. The limit
case of $k=\infty$ was considered in~\cite{DAN,MinNaz,FAA} (where
one can also find more detailed bibliographical remarks). Here we
essentially show that the results obtained there remain valid,
\textsl{mutatis mutandis}, for the case of finite~$k$.

Theorems deriving the asymptotic behavior of $\Ng(x)$ as
$x\to\infty$ from that of $\pig(x)$ are known as (\textit{inverse})
\textit{abstract prime number theorems}, and the corresponding
theorems for $\Ngk(x)$ are a~generalization of these. Apart from
the purely number-theoretic meaning, the function $\Ngk(x)$ has a
natural interpretation in statistical mechanics. Let us enumerate
the elements of $P$ in some way, $P=\{p_1,p_2,\dots\}$, and set
$\la_j=\pa(p_j)$. Then $\Ngk(x)$ is the number of solutions of the
inequality
\begin{equation}\label{01}
\sum_{j=1}^\infty \la_j n_j \le x
\end{equation}
in integers $n_j$ such that
\begin{equation}\label{06}
    0\le n_j < k.
\end{equation}
Inequality~\eqref{01} describes the states with total energy~$\le
x$ of a system of noninteracting indistinguishable particles, $n_j$
being the number of particles at the energy level~$j$ with energy
$\la_j$. Inequalities~\eqref{06} imply that there are at most $k-1$
particles at each energy level. In other words, the particles obey
the Gentile statistics~(see~\cite{Gentile,Gentile1,Khare},
\cite[p.~258]{KvasnikovEquilibrium}), which becomes the well-known
Bose--Einstein statistics (any number of particles at any level)
and Fermi--Dirac statistics (at most one particle at each level) in
the limit cases of $k=\infty$ and $k=2$, respectively. The
logarithm $\ln\Ngk(x)$ is the entropy of the system. Note, however,
that the counting function $\pig(x)$ usually has a power-law
asymptotics rather than the exponential asymptotics~\eqref{e1} in
statistical mechanics, at least if the individual particles have
finitely many degrees of freedom (e.g., see~\cite{Mas08} and the
survey~\cite{Mas13}, where further references can be found). Our
interest in the asymptotics~\eqref{e1} is partly motivated by
problems arising when calculating the number of localized Gaussian
packets in the theory of dynamical systems on metric and decorated
graphs \cite{ChernShaf,ChernTol,Shubert}, where exponential growth
is associated with positivity of topological entropy of the
manifolds in question (see~\cite{Katok,Mane,Pollicott}).

\section{Main results}

Assume that condition~\eqref{e1} is satisfied. The Dirichlet series
\begin{equation}\label{e2}
\zgk=\Sgk \esa, \qquad s=\si+it,
\end{equation}
converges absolutely in the half-plane $\si>1$, and one has the
Euler identity
\begin{equation}\label{e3}
    \zgk=\Pp\frac{1-\espk}{1-\esp}.
\end{equation}
The proof is the same as for the zeta function~$\zg$ of~$G$ (e.g.,
see~\cite[p.~36]{Knopf}), which is the special case of~\eqref{e2}
for $k=\infty$.

Now we are in a position to state the main results of the paper.

\begin{theorem}\label{th1}
Under condition~\eqref{e1}, the function $\Ngk(x)$ has the
following asymptotics as $x\to\infty$\textup:
\begin{equation}\label{e7}
    \Ngk(x)=
    \frac{e^{xs}\zgk}{\sqrt{2\pi(\ln\zgk)''}}\bigg|_{s=\be(x)}
    (1+O(x^{-\ka})),
\end{equation}
where $s=\be(x)>1$ is the unique real solution of the equation
\begin{equation}\label{e5}
    x+(\ln\zgk[s])'=0
\end{equation}
and $\ka>0$ is an arbitrary number such that
\begin{equation}\label{e7a}
     \ka<\frac{\dt}{2+\ga},\qquad\ka\le\frac{1+\ga}{2+\ga}.
\end{equation}
\end{theorem}

Theorem~\ref{th1} gives the asymptotics of $\Ngk(x)$ in terms of
the function $\zgk$, which is itself given by the infinite
product~\eqref{e3}. The formulas for the logarithmic asymptotics
are much simpler and depend only on the constants $\rho$, $\ga$,
and $\dt$ occurring in the asymptotics~\eqref{e1} of the prime
counting function. Namely, the following theorem holds.

\begin{theorem}\label{th2}
Under condition~\eqref{e1}, the function $\ln\Ngk(x)$ has the
following asymptotics as $x\to\infty$\textup:
\begin{equation}\label{entropy-1}
\ln\Ngk(x)=
x+2(\rho\Ga(\ga+2))^{\frac1{\ga+2}}x^{\frac{\ga+1}{\ga+2}}+R(x)
\end{equation}
if $\dt\le\min\{1, 1+\ga\}$, where
\begin{equation*}
    R(x)=O(x^{\frac{\ga+1-\dt}{\ga+2}})\quad
      \text{if $\dt < 1+\ga$ and}
    \quad
    R(x)=O(\ln x)\quad
    \text{if $\dt = 1+\ga$\textup;}
\end{equation*}
\begin{equation}\label{entropy-2}
    \ln\Ngk(x)=
x+2(\rho\Ga(\ga+2))^{\frac1{\ga+2}}x^{\frac{\ga+1}{\ga+2}}
-\frac12\,\frac{\ga+3}{\ga+2}\ln x + O(1)
\end{equation}
if $1\ge \dt > 1+\ga$.
\end{theorem}

\begin{remark}
In contrast to the asymptotics obtained in Theorem~\ref{th1}, the
logarithmic asymptotics provided by Theorem~\ref{th2} does not feel
the difference between the cases of $k=\infty$ and finite~$k$.
\end{remark}

\begin{remark}
Similar results were obtained in \cite{GranovskyStark13} in a
different setting. (The analysis in that paper only applies to the
case in which the mapping $\pa$ is integer-valued.)
\end{remark}

\section{Proof of the theorems}

The proof of both theorems is completely similar to that given
in~\cite{FAA} for the case of $k=\infty$, and here we only give a
brief outline of the reasoning. The argument relies on asymptotic
formulas for the function $\ln\zgk[\si]$ and its derivatives as
$\si \downarrow 1$. These formulas have the form
\begin{align}\label{A11}
\ln\zgk[\si] &= \rho\Ga(1+\ga)(\si-1)^{-1-\ga}(1+o(1)),
\\ \label{A13}
(\ln\zgk[\si])' &= -\rho\Ga(\ga+2)(\si-1)^{-\ga-2}(1+o(1)),
\\ \label{A14}
(\ln\zgk[\si])'' &= \rho\Ga(\ga+3)(\si-1)^{-\ga-3}(1+o(1))
\end{align}
(we only write out the leading terms of the asymptotics) and
coincide modulo $O(1)$ with those obtained in~\cite{FAA} for
$\ln\zg[\si]$, because
\begin{equation*}
    \ln\zgk[\si] -\ln\zg[\si]=
    \sum_{p\in P} \ln(1-e^{-\pa(p) \si k})=O(1)\qquad\text{as
    }\si\downarrow1.
\end{equation*}
(Recall that $k\ge2$.)

\begin{proof}[Outline of proof of Theorem~\ref{th1}]
We have
\begin{equation}\label{def_N}
\Ngk(x) = \Sgk  H\bigg(\frac{x - \pa(a)}{\e} \bigg),
\end{equation}
where $H(x)$ is the Heaviside step function and $\e>0$ is
arbitrary. Take smooth functions $\chi_\pm(x)$ such that
\begin{equation}\label{chipm}
    \chi_-(x)\le H(x)\le\chi_+(x)\quad\text{for all $x$,}\quad
    \chi_\pm(x)=H(x),\quad \abs{x}\ge1;
\end{equation}
then
\begin{equation*}
    \Sgk\chi_-\BL(\frac{x-\pa(a)}{\e}\BR)\le
    \Ngk(x)\le \Sgk\chi_+\BL(\frac{x-\pa(a)}{\e}\BR).
\end{equation*}
Using the generalization in~\cite[Proposition~3]{FAA} of the
well-known Perron formula~\cite[p.~12, Theorem~13]{HR}, we obtain
\begin{equation}\label{e13}
    I_-(x,\e)\le
    \Ngk(x)\le
    I_+(x,\e),
\end{equation}
where
\begin{equation}\label{e14a}
    I_\pm(x,\e) =
    \frac1{2\pi i}\int_{\si-i\infty}^{\si+i\infty}
    e^{xs}\zgk\; \e\wt\chi_\pm(\e s)\,ds,
\end{equation}
$\wt\chi_\pm(s)$ is the two-sided Laplace transform of
$\chi_\pm(x)$, and the integrals are independent of the choice of
$\si>1$.

Now we compute these integrals by the mountain pass
method~\cite[Ch.~4, p.~170]{Fed}. The phase function is
\begin{equation}\label{phase0}
    S(x,s)=xs+\ln\zgk,
\end{equation}
and the amplitude is $\e \wt\chi_\pm(\e s)$. The equation $\pa
S(x,s)/\pa s=0$ for the stationary points of the phase
function~\eqref{phase0} coincides with~\eqref{e5} and has a unique
real solution $s=\be(x)>1$ for each $x>0$. Further, $\be(x)\to1$ as
$x\to\infty$, and it follows from \eqref{A13} that
\begin{equation}\label{e6}
  \be(x) - 1 \sim C x^{-1/(\ga+2)}, \qquad C = \big(\rho
   \Gamma(\ga+2)\big)^{\frac1{\ga+2}}.
\end{equation}
In~\eqref{e14a}, we take the integration contour to be given by
$\si = \be(x)$; this is a mountain pass contour for this integral.
We make a change of the integration variable $s$ by the formula $s
= \be(x)+i\xi$, so that the contour of integration with respect to
the new variable~$\xi$ coincides with the real line. Further, set
\begin{equation}\label{zamena}
    x=x(\be)\equiv-(\ln\zgk[\be])';
\end{equation}
this is the inverse function of~$\be(x)$. These transforms give the
integrals
\begin{align}\label{newint}
    I_\pm(x(\be),\e)&= \frac{1}{2\pi}
    \int_{-\infty}^\infty
    e^{S(x(\be),\be+i\xi)} \;
    \e\wt\chi_\pm(\e(\be+i\xi)) \,d\xi,
\end{align}
with the saddle point $\xi = 0$. By using the Mountain Pass
Theorem~\cite[Ch.~4, Theorem~1.3, p.~170]{Fed} and the estimates
in~\cite{FAA} for the phase function $\Phi(\be,\xi) =
S(x(\be),\be+i\xi)$ and the amplitude $\ph_\pm(\be,\xi) = \e(\be)
\wt\chi_\pm(\e(\be) (\be+i\xi))$ with $\e(\be) =
(\be-1)^{\ka(2+\ga)},$ $\ka>0$, we finally obtain
\begin{equation}\label{fine}
\begin{split}
    I_\pm(x(\be),\e(\be))&=
    \frac{e^{x(\be)\be}\zgk[\be]}{\sqrt{2\pi (\ln\zgk[\be])''}} \\
    &\qquad{}\times
    (1+O(\e(\be))+O(\be-1)+O((\be-1)^{1+\ga})),
\end{split}
\end{equation}
which, together with \eqref{e13}, gives \eqref{e7}.
\end{proof}

\begin{proof}[Outline of proof of Theorem~\ref{th2}]
This theorem follows if one takes the logarithm of both sides of
formula~\eqref{e7} in Theorem~\ref{th1} and then uses the
asymptotics~\eqref{e6} of $\be(x)$ together with the asymptotics of
$\ln\zgk[\si]$ and $(\ln\zgk[\si])''$ whose leading terms are given
by~\eqref{A11} and~\eqref{A14}, respectively.
\end{proof}

\section{Simulation results}\label{s4}

\begin{figure}\centering
\includegraphics[width=150pt,angle=0]{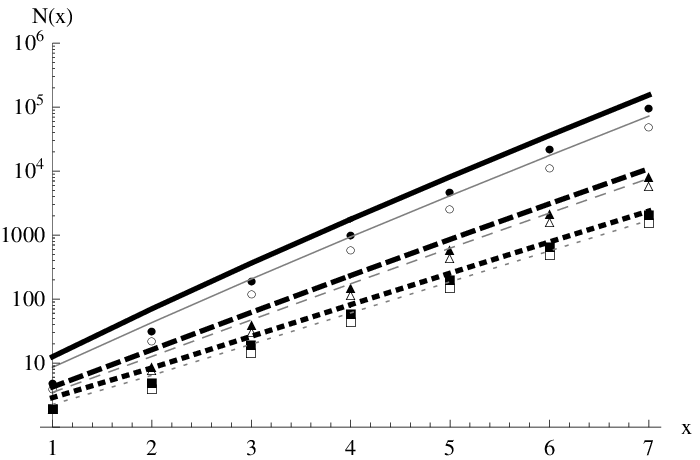}
\includegraphics[width=150pt,angle=0]{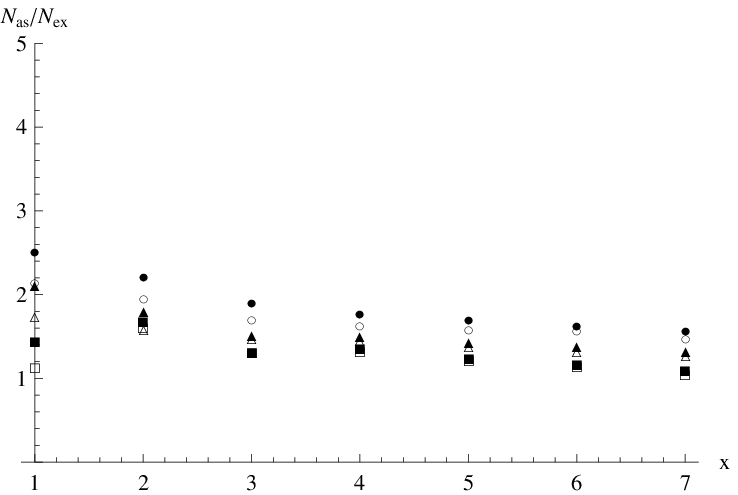}
\includegraphics[width=150pt,angle=0]{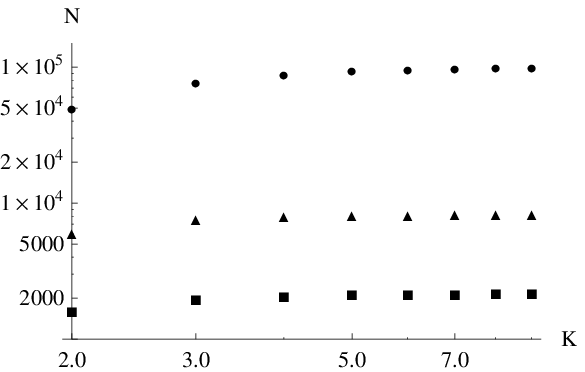}
\caption{Left: the comparison of exact (points) and asymptotic
(lines) values of $\Ngk(x)$ for $k = \infty$ (solid points and
black lines) and $k = 2$ (empty points and gray lines). The
parameter values are $\rho = 0.5, 1, 2$ (dotted lines and squares,
dashed lines and triangles, and solid lines and circles,
respectively).
\newline
Right: the ratio of the asymptotic values to the exact values of
$\Ngk$ for $k=2, \infty$ and $\rho = 0.5, 1, 2$ with the same
notation.
\newline
Bottom: the dependence of $\Ngk(x)$ on the parameter $k\in[2,8]$ at
the point $x=7$ for $\rho = 0.5, 1, 2$ (squares, triangles, and
circles, respectively).} \label{Fig_N(x)}
\end{figure}

To illustrate the results, consider the arithmetical semigroup $G$
with primes $p_n\in P$, $n=1,2,\dots$, and with $\pa(p_n) =
\ln(\frac {n + \rho}{\rho})$. Then
\begin{equation*}
   \pig(x) = [\rho e^x - \rho] = \rho e^x + R(x), \qquad
  -\rho \le R(x) \le -\rho + 1.
\end{equation*}
We compare the exact values of $\Ngk(x)$ and the asymptotic values
given by Theorem~\ref{th1} for $\rho = 0.5, 1$, and $2$ at the
points $x=1,2,...,7$. The results are presented in
Fig.~\ref{Fig_N(x)} for $k=2$ and $k=\infty$ (the Fermi and Bose
cases). We also present the dependence of $\Ngk$ on the parameter
$k \in [2,8]$ at the point $x=7$.

The figure illustrates the convergence of asymptotic formulas to
the exact values (right). It also illustrates the fact, that the
rate of growth $\Ngk(x)$ with $x$ is the same for all $k$ and
asymptotics $\Ngk(x)$ differs by a factor (left). This fact follows
from corollaries. We also can see from the bottom figure that this
factor tends to $1$ very rapidly with increasing parameter $k$ and
the value $k=10$ can already be treated as ``infinity'' from the
viewpoint of convergence.

\subsection*{Acknowledgements}

The research of the first author was carried out in the framework
of the Academic Fund Program at the National Research University
Higher School of Economics (HSE) in 2015--2016 (grant
no.~15-01-0091) and
supported within the framework of a subsidy granted to the HSE
by the Government of the Russian Federation for the implementation of the
Global Competitiveness Program.

The authors are grateful to V.\,P.~Maslov and S.\,Yu.~Dobrokhotov for
encouragement and also to A.\,I.~Shafarevich, A.~Tolchennikov, and
D.~Yakobson for very useful discussions.

\end{document}